\documentstyle[11pt]{article}
\input amssymb.sty
\newcommand{\ti}[1]{\tilde{#1}}

\newcommand{\Om}{\Omega}
\newcommand{\de}{\delta}
\newcommand{\al}{\alpha}
\newcommand{\te}{\theta}

\newcommand{\be}{\beta}

\newcommand{\D}{\Delta}
\newcommand{\ve}{\varepsilon}

\newcommand{\vf}{\varphi}
\newcommand{\G}{\Gamma}
\newcommand{\ka}{\kappa}

\newcommand{\ga}{\gamma}

\newcommand{\mat}[4]{\left(\begin{array}{cc}{#1}&{#2}\\{#3}&{#4}
\end{array}\right)}
\newcommand{\beq}[1]{\begin{equation}\label{#1}}
\newcommand{\eq}{\end{equation}}
\newcommand{\beqn}[1]{\begin{eqnarray}\label{#1}}
\newcommand{\eqn}{\end{eqnarray}}
\newcommand{\p}{\partial}
\newcommand{\di}{{\rm diag}}
\newcommand{\oh}{\frac{1}{2}}

\newcommand{\sL}{{\rm sl}(2,{\mathbb C})}
\newcommand{\SL}{{\rm SL}(2,{\mathbb C})}
\def\SU{{\rm SU}(2)}

\def\f1#1{\frac{1}{#1}}
\newcommand{\rar}{\rightarrow}

\newcommand{\bz}{\bar{z}}
\newcommand{\bq}{\bar{q}}

\newcommand{\bae}{\bar{e}}
\newcommand{\bh}{\bar{h}}
\newcommand{\baf}{\bar{f}}
\newcommand{\baka}{\bar{\kappa}}
\newcommand{\ot}{\otimes}

\newtheorem{predl}{Proposition}[section]

\newtheorem{cor}{Corollary}[section]

\begin{document}
\vspace{0.3in}
\begin{flushright}
ITEP-TH-57/01
\end{flushright}
\vspace{10mm}
\begin{center}
{\Large {\bf Unitary Representations of Quantum Lorentz Group\\
and Quantum Relativistic Toda Chain}}\\
\vspace{5mm}
M.A.Olshanetsky \\
{\sf ITEP, 117259, Moscow, Russia; MPIM, Vivatgasse 7, 53111, Bonn, Germany}\\
{\sl e-mail olshanet@mpim-bonn.mpg.de}\\
V.-B.K.Rogov \\
{\sf MIIT, 127994, Moscow, Russia} \\
{\sl e-mail vrogov@cemi.rssi.ru}\\
\vspace{5mm}
2001\\
\end{center}
\begin{abstract}
The aim of this paper is to give a group theoretical interpretation of
the three types of Bessel-Jackson functions. We consider a family of
quantum Lorentz groups and a family of quantum Lobachevsky spaces.
For three members of quantum Lobachevsky spaces the Casimir operators
give rise to the two-body relativistic open Toda lattice Hamiltonians.
Their eigen-functions are the modified Bessel-Jackson functions of three
types. We construct the principal series of unitary irreducible representations
of the quantum Lorentz groups. Special matrix elements in the irreducible spaces
are the Bessel-Macdonald-Jackson functions. They are the wave functions
of the two-body relativistic open Toda lattice. We obtain integral
representations for these functions.
\end{abstract}

\section{Introduction}
There exist deep interrelations between the Harmonic Analysis on symmetric spaces
and quantum integrable systems \cite{OP}.
In this approach the zonal spherical functions play the role of the wave-functions
of some integrable models. In particular, the Bessel functions are the zonal spherical
functions related to the group symmetries of the Euclidean spaces \cite{Vi}.
On the other hand, they are the wave-functions of the rational Calogero-Moser
model.

The modified Bessel functions and the Bessel-Macdonald functions arise in a
different construction. Consider the eigen-functions of the Laplace-Beltrami operator
in the horospheric coordinates on the Lobachevsky space. The modified Bessel
functions and the Bessel-Macdonald functions are their Fourier transform on
a horosphere \cite{Vi}. From the point of view of integrable systems  the Laplace-Beltrami
operator in horospheric coordinates gives rise to the simplest form of the open
Toda Hamiltonian.
B. Kostant generalized this connection and established the similar relations
between an arbitrary open Toda lattice and the Whittaker model of irreducible
representations of the splitted simple group $G$ \cite{K} (see also reviews
\cite{STS,GMM}).
This construction can be generalized in different ways. For example,
the periodic Toda lattices arise if one replaces a simple group $G$ by
the central extended loop group \cite{FF,GW}.

The zonal spherical functions for the quantum groups were investigated in
numerous papers (see \cite{NYM,V}). In particular,
the $q$-Bessel-Jackson functions were
described by Vaksman and Korogodsky \cite{VK} as
the zonal spherical functions related to the quantum plane.

In our paper \cite{OR1}  we applied the Kostant scheme to
 the quantum Lorentz group constructed by Podlez and Woronowicz \cite{PW}.
It turns out that the corresponding integrable system is the open two-body
relativistic Toda lattice \cite{Ru}. Its Hamiltonian is related to the Casimir
operator of  the quantum Lorentz group being written in an analog of the
horospheric coordinates on the corresponding quantum Lobachevsky space.
Up to the $q$-exponents multipliers the eigen-functions are the modified
$q$-Bessel-Jackson and the $q$-Bessel-Jackson-Macdonald functions. Their
properties were investigated in \cite{OR2,OR3,R1}. For general quantum groups
an analog of the Whittaker model was constructed in \cite{Et,Se}. Recently, the
integral representations for the wave functions of the $N$-body case based on the
Whittaker model  of $U_q({\rm sl}_N({\mathbb R}))$ was presented in \cite{KhL}.

Here we come back to the two-body case and investigate it in detail.
Our analysis is based on the Whittaker model for the
quantum Lorentz group $U_q({\rm sl}_2({\mathbb C}))$.
In this way we describe the group theoretical approach to the modified $q$-Bessel-Jackson
and the $q$-Bessel-Jackson-Macdonald functions.
For this purpose we construct unitary irreducible representations of
a twisted family of quantum Lorentz groups.
As by product we obtain some results in the Harmonic Analysis on quantum
Lobachevsky spaces.
The quantum Lobachevsky spaces play the role of homogeneous spaces with respect to the
twisted quantum Lorentz groups. Three members from the family
are distinguished. For these cases the Casimir operators
realized on the quantum Lobachevsky spaces lead to second order difference
operators. Up to a conjugation they coincide with the quantum relativistic
two-body open Toda Hamiltonians. Their eigen-functions are the three
types of the modified Bessel-Jackson functions. We construct an analog of the
principle series unitary irreducible representations.
Then we consider some special matrix elements of group elements acting in an
irreducible space of the class-one representations.
They are the wave functions of the Toda lattice. These matrix elements have the
form of the double integral. It follows from \cite{OR3} that they are
$q$-Bessel-Jackson-Macdonald functions.

The paper has the following structure. In Section 2 we recall interrelations between
the Lorentz group $\SL$ and the two-body open Toda model. In Section 3 we
construct the twisted family of quantum Lorentz groups and the corresponding family
of quantum Lobachevsky spaces. The principle series of unitary irreducible
representations are constructed in Section 4. The eigen-functions of the Casimir operators
on the quantum Lobachevsky spaces are constructed in Section 5.
The integral representations of the wave functions for the open Toda model
are presented in Section 6. Finally, in Section 7 we give the representation for the
wave-functions as the Mellin-Barns integral. The similar representation was
obtained in \cite{KhL}.

{\bf Acknowledgments}\\
{\sl  The work of M.O. is supported in part by 
INTAS-99-01782, RFBR-00-01-00143 and RFBR-00-15-96455 for the support of
scientific schools.
The work of V.R. is supported in part by 
RFBR-00-01-00143. M.O. is 
 grateful to the Max-Planck-Institut f\"{u}r Mathemamatik in Bonn
for the hospitality, where this paper was completed.
The authors thank S.Kharchev for valuable discussions.

}

\vspace{10mm}

\section{Classical Case}
\setcounter{equation}{0}

{\bf 1. Lorentz group.}
The matrix Lorentz group is the group of complex
unimodular matrices of second order
$$
G=\SL=\mat{\al}{\be}{\ga}{\de},~(\al\de-\be\ga=1).
$$
Its Lie algebra  has  three generators
$$
e=\mat{0}{1}{0}{0},~h=\oh\di(1,-1),~f=e^T
$$
with commutators
$$
[e,f]=2h,~[h,e]=e,~[h,f]=-f.
$$
They generate ${\rm sl}(2,{\mathbb R})$ over ${\mathbb R}$ and $\sL$ over ${\mathbb C}$.
As a real algebra $\sL$ has six generators and one should consider the second
 copy of  ${\rm sl}(2,{\mathbb R})$ generated by $(\bae,\bar{h},\baf)$ .

\bigskip
{\bf 2 Principal series} \cite{GGV}.
The principal series $\pi_{\nu,n}$ of the unitary representations of $G=\SL$
are defined in the following way.

Let $N^-$ be the subgroup of the lower triangular matrices,
$A$ are the real diagonal matrices and $M$ are the unitary diagonal matrices.
The Borel subgroup $B$ of $G$
$$
B=N^-AM,~~b=v\cdot\di(h,h^{-1})\di(e^{2\pi i\te},e^{-2\pi i\te}),~~(v\in N^-)
$$
has characters
$$
\chi_{\nu,n}(vam)=\exp\{(i\nu-1)\log h+n\te\}.
$$
The principal series representations $\pi_{\nu,n}$ of $G$ are induced
by the characters of $B$ in the space of smooth functions on $G$
\beq{2.1}
f(bg)=\chi_{\nu,n}(b)f(g), ~~\pi_{\nu,n}(g)f(x)=f(xg),~x\in\SL.
\eq
It is amount to look on the action in the space of sections
of a linear bundle over ${\bf P}^1$. To this end
consider the Gauss decomposition of the dense subset
$\{g|\al\neq 0\}\subset\SL$
$$
G=BN,~~(N -{\rm nilpotent~subgroup~of~the~ upper~ triangular~matrices}).
$$
\beq{2.2}
g=bn(g),~n_{12}(g)=\frac{\be}{\al}=z,~~b_{11}=\al.
\eq
Then (\ref{2.2}) means that the right action (\ref{2.1}) of $\SL$  induces
the M\"{o}bius transform on ${\bf P}^1$
$$
z\rar zg=\frac
{\de z+\be}{\ga z+\al}.
$$
Consider the linear bundle
${\cal L}_{\nu,n}$ over ${\bf P}^1$ with the space of sections
\beq{2.3}
V_{\nu,n}=\G({\cal L}_{\nu,n})\sim{\cal A}^{(-\frac{r}{2},-\frac{\ti{r}}{2})}({\bf P}^1)
~(r=i\nu+n-1,~\ti{r}=i\nu-n-1).
\eq
The sections have the form
$f(z,\bz)(dz)^{-\frac{r}{2}}(d\bz)^{-\frac{\ti r}{2}}$, where $f(z,\bz)$
are smooth functions with the asymptotic
\beq{2.4}
f(z,\bz)_{z\to\infty}\sim z^r\bz^{\ti{r}}.
\eq
By means of (\ref{2.2}) we can reduce $\pi_{\nu,n}(g)$ (\ref{2.1}) to the
action on $V_{\nu,n}$ (\ref{2.3})
\beq{2.5}
\pi_{\nu,n}(g)f(z,\bz)=(\ga z+\al)^{i\nu+n-1}
(\bar{\ga} \bz+\bar{\al})^{i\nu-n-1}
f(zg,\bz\bar g).
\eq

Since $\bar{r}+\ti{r}=-2$ there is a Hermitian form on
${\cal A}^{(-\frac{r}{2},-\frac{\ti{r}}{2})}({\bf P}^1)$
\beq{2.6}
<f_1|f_2>
=\int_{{\bf P}^1}f_1(z,\bar{z})\overline{f_2(z,\bar{z})}dzd\bar{z}.
\eq

The representation $\pi_{\nu,n}$ (\ref{2.5})  is realized in $V_{\nu,n}$
by the unitary operators because
$$
<\pi_{\nu,n}(g)f_1|f_2>=<f_1|\pi_{\nu,n}(g^{-1})f_2>.
$$

The infinitesimal version of this construction takes the form:
$$
e\rar T^+=\p_z,~~h\rar T^3=-z\p_z+\frac{r}{2},
$$
\beq{2.7}
~f\rar T^-=-z^2\p_z+rz,~~(r=n-1+i\nu)
\eq
$$
\bar{e}\rar \bar{T}^+=\p_{\bar{z}},~
\bar{h}\rar \bar{T}^3=-\bar{z}\p_{\bar{z}}+\frac{\ti{r}}{2},
$$
\beq{2.8}
\bar{f}\rar \bar{T}^-=-\bar{z}^2\p_{\bar{z}}+\ti{r}\bar{z},~~
(\ti{r}=-n-1+i\nu).
\eq
This action preserves the asymptotic (\ref{2.4}).

There are two Casimir operators of $\SL$
$$
\Om=h^2+h+fe, ~~\ti{\Om}=\bar{h}^2+\bar{h}+\bar{f}\bar{e}.
$$
They become the scalar operators in $V_{\nu,n}$
\beq{2.9}
\Om|\Psi>=(\f1{4}(i\nu+n)^2-\f1{4})|\Psi>,~~
\ti{\Om}|\Psi>=(\f1{4}(i\nu-n)^2-\f1{4})|\Psi>,~~
|\Psi>\in V_{\nu,n}.
\eq

If $n=0$ then there exists a $\SU$ invariant vector  in $V_{\nu,0}$
$$
(1+|z|^2)^{i\nu-1}.
$$

\bigskip
\noindent
{\bf 3. Lobachevsky space.}
A subclass of these representations $r=i\nu-1,~(n=0)$
(the representations of class one) are related to the decomposition
of functions on the Lobachevsky space
 $ {\bf L}=\SU\backslash \SL$. The later can be realized
as the space of second
order unimodular Hermitian positive definite matrices.
Any $x\in {\bf L}$ can be represented as
\beq{2.10}
x=g^{\dagger}g=\mat{\bar{\al}\al+\bar{\ga}\ga}{\bar{\al}\be+\bar{\ga}\de}
{\bar{\be}\al+\bar{\de}\ga}{ \bar{\be}\be+\bar{\de}\de}.
\eq

The Iwasawa decomposition
$$
g=kb, ~g\in  \SL,~ k\in \SU,~b\in AN
$$
allows to introduce the horospheric coordinates on ${\bf L}$.
If
$$b=\mat{h}{hz}{0}{h^{-1}},$$
then from (\ref{2.10})
\beq{2.11}
x=b^{\dagger}b=
\mat{\bar{h}h}{\bar{h}hz}{\bar{z}\bar{h}h}{ \bar{z}\bar{h}hz+(\bar{h}h)^{-1}}.
\eq
The triple $(H=\bar{h}h,z,\bar{z})$ is uniquely determined by $x$. It is called
the horospheric coordinates of  $x$.
 It follows from (\ref{2.10}) and (\ref{2.11}) that
$$H=\bar{\al}\al+\bar{\ga}\ga,$$
$$Hz=\bar{\al}\be+\bar{\ga}\de,$$
$$\bar{z}H=\bar{\be}\al+\bar{\de}\ga.$$

Consider  the space of smooth complex valued integrable functions
${\cal R}_{\bf L}=\{f(H,z,\bar{z})\}$ on ${\bf L}$:
\beq{2.12}
\int |f(H,z,\bar{z})|^2HdHdzd\bar{z}<\infty.
\eq
The  Lie operators generated by the right shifts act on
${\cal R}_{\bf L}$ as
$$
e\rar {\cal D}^+=\partial_z,
$$
\beq{2.13}
h\rar {\cal D}^3=\frac{1}{2}H\partial_H-z\partial_z,
\eq
$$
f\rar {\cal D}^-=Hz\partial_H-z^2\partial_z+H^{-2}\partial_{\bar{z}}.
$$
The generators $\bae, \bh, \baf$ give rise to the complex conjugate operators
$$
\bae\rar  \bar{\cal D}^+=\p_{\bz},
$$
\beq{2.14}
\bh\rar \bar{\cal D}^3=\frac{1}{2}H\partial_H-\bz\partial_{\bz},
\eq
$$
\baf\rar \bar{\cal D}^-=H\bz\partial_H-\bz^2\partial_{\bz}+
H^{-2}\partial_{\bz}.
$$

Consider a family $V_\nu, ~\nu\in{\mathbb R}$ of  subspaces
in ${\cal R}_{\bf L}$
$$
V_\nu=\{f(\bz,aH,z)=a^{(i\nu-1)}f(\bz,H,z)\}.
$$
Assume that $e,h,f$ act on the subspace of holomorphic ($\bz$ independent)
functions and $\bae,\bh,\baf$ act
on the subspace of antiholomorphic functions.
After comparison (\ref{2.13}) and (\ref{2.14}) with (\ref{2.7}) and
(\ref{2.8}) one concludes that $V_\nu=V_{\nu,0}$.
 The invariant integral (\ref{2.12}) coincides with
the hermitian form (\ref{2.6}) on the irreducible subspace.

The Casimir operators in the horospheric coordinates take the form
\beq{2.15}
\Omega=\frac{1}{4}H^2\partial^2_H+\frac{3}{4}H\partial_H+
H^{-2}\partial^2_{z\bar{z}},
\eq
$$
\ti{\Om}=\Om.
$$
According to (\ref{2.9}) they are scalar operators on $V_\nu$
$$
\Om|_{V_\nu}=-\f1{4}(\nu^2+1){\rm Id}.
$$

\bigskip
\noindent

{\bf 4.Whittaker functions}.
Let $U_{\nu,\mu}=\{f_\nu(\bz,H,z)\}$  be the space of smooth functions on $ {\bf L}$
satisfying the following conditions:\\
(i) $f_\nu(\bz,H,z)$ are the eigen-functions of the Casimir operator (\ref{2.15})
$$
\Om f_\nu(\bz,H,z)=-\f1{4}(\nu^2+1)f_\nu(\bz,H,z).
$$
(ii) $\p_z f(\bz,H,z)=i\mu f(\bz,H,z),~~
\p_{\bz} f(\bz,H,z)=i\mu f(\bz,H,z)$.

The last condition means that
$$
f_\nu(\bz,H,z)=\exp i\mu (z+\bz)F_\nu(H).
$$

Substituting this expression in the equation (i)
one obtains
\beq{2.16}
(\frac{1}{4}H^2\partial^2_H+\frac{3}{4}H\partial_H-
H^{-2}\mu^2)F_\nu(H)= -\f1{4}(\nu^2+1)F_\nu(H).
\eq
Two independent solutions of this equation are modified Bessel functions
$$H^{-1}I_{\pm i\nu}(2\mu H^{-1}).$$
The bounded solution is the Bessel-Macdonald function
\beq{2.17}
H^{-1}K_{i\nu}(2\mu H^{-1})=
H^{-1}\frac{\pi(I_{-i\nu}(2\mu H^{-1})-I_{i\nu}(2\mu H^{-1}))}
{2i\sinh\nu}.
\eq
The last expression allows to find the asymptotic for $H\rar 0$
$$
K_{i\nu}(2\mu H^{-1})\sim \frac{\pi}{2i\sinh \nu}
\left(\f1{\G(1-i\nu)}\left(\frac{H}{2}\right)^{-i\nu}-
\f1{\G(1+i\nu)}\left(\frac{H}{2}\right)^{i\nu}\right).
$$

The bounded solution (\ref{2.17}) can be represented as a special matrix element
of the operator
\beq{2.18}
g(T^3|\vf)=\exp 2(T^3\otimes \vf)\in{\cal U}_\SL\otimes{\cal A}_{\bf L},
~~~H=e^\vf,
\eq
which is a special group element  defined in the tensor product of
the universal enveloping algebra ${\cal U}_\SL$ and the group algebra
${\cal A}_{\bf L}$. It acts on the irreducible space $V_\nu$.
Define two vectors $\psi_L,~~\psi_R$ such that
$$
T^+\psi_R=i\mu \psi_R, ~~~\bar{T}^+\psi_R=i\mu \psi_R,
$$
$$
 (T^+-\bar{T}^-)\psi_L =(T^3-\bar{T}^3)\psi_L=0.
$$
It means that $\psi_R$ is the eigen-vector of the  nilpotent subgroup $N$,
and $\psi_L$ is the $SU(2)$-invariant vector. The state $\psi_R$ is called
{\sl the Whittaker vector}. Then it follows from (\ref{2.7}), (\ref{2.8}) that
\beq{2.19}
\psi_R=\exp i\mu(z+\bz),
\eq
 and
\beq{2.20}
\psi_L=\f1{(1+|z|^2)^{ 1 -i\nu }}.
\eq
In fact, $\psi_R$ does not lie in $V_\nu$, but rather it is a distribution over some
subspace of $V_\nu$. Nevertheless, the integral
\beq{2.21}
<\psi_L|g(T^3|\vf)|\psi_R>
\eq
is well defined. This construction can be generalized to an arbitrary
splitted Lie group \cite{K,STS}. The bounded solution is called {\sl the
Whittaker function}.

It will be instructive to prove that this matrix element  satisfies (\ref{2.16}).
In fact, taking in account that
$g(T^3|\vf)T^-=e^{2\vf}T^- g(T^3|\vf)$ one has
$$
-\f1{4}(\nu^2+1)<\psi_L|g(T^3|\vf)|\psi_R>=
<\psi_L|g(T^3|\vf)\Om|_{V_\nu}|\psi_R>=
$$
$$
(\frac{1}{4}H^2\partial^2_H+\frac{3}{4}H\partial_H) <\psi_L|g(T^3|\vf)|\psi_R>
+<\psi_L|g(T^3|\vf)T^-T^+|\psi_R>.
$$
Then the conditions (\ref{2.19}), (\ref{2.20}) and the commutation relation
$H\p_H g(T^3|\vf)= g(T^3|\vf)T^3$ give
$$
<\psi_L|g(T^3|\vf)T^-T^+|\psi_R>=-\mu^2H^{-2}<\psi_L|g(T^3|\vf)|\psi_R>.
$$
Thus, we reproduce the left hand side of (\ref{2.16}).

Substituting the explicit expressions (\ref{2.18}),(\ref{2.19}),(\ref{2.20})
to (\ref{2.21}) one finds the bounded solution of the equation (\ref{2.16})
\beq{2.22}
<\psi_L|g(T^3|\vf)|\psi_R>=
\exp((i\nu+1)\vf)\int\int d\bz dz
\frac{\exp(i\mu(z+\bz))}
{(1+e^{2\vf}|z|^2)^{1+i\nu}}\sim H^{-1} K_{i\nu}(2\mu H^{-1}).
\eq
In this formula the
angular integration in the polar coordinates  allows to rewrite it as
\beq{2.23}
H^{-1} K_{i\nu}(2\mu H^{-1})=
\G(1+i\nu) H^{-i\nu-1}\int_0^\infty\frac
{J_0(2\mu\rho)\rho}
{(H^{-2}+\rho^2)^{1+i\nu}}d\rho.
\eq

Let us substitute $f_\nu(H)=\ e^{-\vf}\Phi_\nu(\vf)$ in (\ref{2.16}).
Then (\ref{2.16}) takes the form of the quantum Liouville equation
\beq{2.24}
\left(\frac{1}{2}\p^2_\vf-2\mu^2e^{-2\vf}\right)\Phi(\vf)=-\frac{1}{2}\nu^2\Phi(\vf).
\eq

The equation (\ref{2.16}) has a special simple form in the momentum representation.
Assume that $Hf(H)=\Phi(e^{\vf})$ and take the Fourier transform
$$
\psi_\nu(p)=\f1{2\pi i}\int_{-\infty}^{\infty}\Phi(\vf)e^{ip\vf}d\vf.
$$
The equation (\ref{2.16}) for $\psi_\nu(p)$ takes the form
$$
\frac{1}{4}(\nu^2-p^2)\psi_\nu(p)-\mu^2\psi_\nu(p-2i)=0.
$$
The bounded solution to this equation vanishing when $p\rar\pm\infty$ is
\beq{2.25}
\psi_\nu(p)=a(\nu)\mu^{ip}\G(\frac{i}{2}(p+\nu))\G(\frac{i}{2}(p-\nu)).
\eq
\vspace{10mm}

\section{Quantum Lorentz group}
\setcounter{equation}{0}

{\bf 1. General construction.}
Deformations of the group algebra ${\cal A}(\SL)$ based on the
Iwasawa and the Gauss decompositions were
considered in \cite{PW} and \cite{WZ1}. Later three parameter
deformations were found in \cite{WZ2}.
We consider here the dual object ${\cal U}_q(\SL)$ and describe a two parameter
family ${\cal U}^{(r,s)}_q(\SL)$.

We start from a pair of the standard ${\cal U}_q({\rm SL}_2)$ Hopf algebra.
The first one is generated by $A,B,C,D$ and the unit with
the relations
$$
AD=DA=1,~AB=qBA,~BD=qDB,
$$
\beq{3.1}
AC=q^{-1}CA,~CD=q^{-1}DC,
\eq
$$
[B,C]=\frac{1}{q-q^{-1}}(A^2-D^2).
$$
It is a Hopf algebra where the coproduct is
defined as
$$\D(A)=A\otimes A,~\D(D)=D\otimes D,$$
\beq{3.2}
\D(B)=A\otimes B+B\otimes D,
\eq
$$
\D(C)=A\otimes C+C\otimes D,
$$
with the counit
$$\ve\mat{A}{B}{C}{D}=\mat{1}{0}{0}{1},$$
and the antipode
$$
S\mat{A}{B}{C}{D}=\mat{D}{-q^{-1}B}{-qC}{A}.
$$
There is a copy of this algebra ${\cal U}^*_q({\rm SL}_2)$ generated by  $A^*,B^*,C^*,D^*$
with the relations following from (\ref{3.1}). The star generators commute with $A,B,C,D$.
The coproduct for the star generators is determined by $\D(X^*)=(\D(X))^*$
and the antipode is extracted from the rule $S\circ*\circ S\circ*=id$.
The pair of two independent copies ${\cal U}_q({\rm SL}_2)$ and ${\cal U}^*_q({\rm SL}_2)$,
which  commute and cocommute,
defines the untwisted version of the quantum Lorentz algebra ${\cal U}_q(\SL)$.
Note that this form   differs  from
the Iwasawa \cite{PW} and the Gauss \cite{WZ1}
constructions.

The ${\cal U}_q(\SU)$ subalgebra \cite{FRT} is defined as
\beq{3.3}
 {\cal U}_q(SU_2)=\{B^*=C,~~C^*=B,~~A^*=A\}.
\eq

There are two Casimir elements in ${\cal U}^{(r,s)}_q(\SL)$ which commute with
any $u\in{\cal U}_q^{r,s}(\SL)$.
\beq{3.6}
\Omega_q:=\frac{(q^{-1}+q)(A^2+A^{-2})-4}{2(q^{-1}-q)^2}+
\frac{1}{2}(BC+CB)
\eq
\beq{3.7}
\ti{\Omega}_q:=\frac{(\bq^{-1}+\bq)(A^{*2}+A^{*-2})-4}{2(\bq^{-1}-\bq)^2}+
\frac{1}{2}(B^*C^*+C^*B^*)
\eq

We use the Casimir operator in the equivalent form
\beq{3.8}
\Om_q=\frac{qA^2+q^{-1}A^{-2}-2}{(q-q^{-1})^2}+CB.
\eq

To go further we following V.Drinfeld \cite{D}
introduce depending on two parameters  twist in the form \cite{Re}.
Let $A=q^K$, $A^*=q^{K^*}$ and
$$
{\cal F}(r,s)=q^{rK^*\ot K+K\ot sK^*}\in {\cal U}_q({\rm SL}_2)\ot
{\cal U}^*_q({\rm SL}_2).
$$
The element ${\cal F}{(r,s)}$ is the so-called left twisting because
it satisfies
$$
({\cal F}\ot 1)(\D\ot id){\cal F}=(1\ot{\cal F})(id\ot\D){\cal F},
$$
$$
(\ve\ot id){\cal F}=1=(id\ot\ve){\cal F}.
$$
These properties allow to generate a new coproduct from the old one:
$$
\D^{\cal F}={\cal F}\D{\cal F}^{-1}.
$$
In our case
$$
\D^{\cal F}(A)=\D(A)=A\ot A,
$$
\beq{3.4}
\D^{\cal F}(B)=(A^*)^{-r}A\otimes B+B\otimes D(A^*)^s,
\eq
$$
\D^{\cal F}(C)=(A^*)^rA\otimes C+C\otimes D(A^*)^{-s}.
$$
The antipode is changed in a consistent way with the new
coproduct.
$$
m(S^{\cal F}\otimes id)\D^{\cal F}(u)=
m(id\otimes S^{\cal F})\D^{\cal F}(u)=\ve(u),
$$
where $m$ is the multiplication. So we have
$$
S^{\cal F}\mat{A}{B}{C}{D}=
\mat{D}{-q^{-1}(A^*)^{r-s}B}{-q(A^*)^{s-r}C}{A}.
$$
In this way we obtain the twisted algebra ${\cal U}^{\cal F}_q({\rm SL}_2)$.
 The algebra $({\cal U}_q^*)^{\cal F}({\rm SL}_2)$
with ${\cal F}$ depending on the same pair
$(r,s)$ can be constructed in the similar way.
$$
\D^{\cal F}(A^*)=\D(A^*)=A^*\ot A^*,
$$
\beq{3.5}
\D^{\cal F}(B^*)=A^{-r}A^*\otimes B^*+B^*\otimes D^*A^s,
\eq
$$
\D^{\cal F}(C^*)=A^rA^*\otimes C^*+C^*\otimes D^*A^{-s}.
$$
$$
S^{\cal F}\mat{A^*}{B^*}{C^*}{D^*}=
\mat{D^*}{-\bq A^{r-s}B^*}{-\bq^{-1}A^{s-r}C^*}{A^*}.
$$

The pair ${\cal U}^{\cal F}_q({\rm SL}_2)$ and $({\cal U}_q^*)^{\cal F}({\rm SL}_2)$
forms the family of commutative, noncocommutative
quantum Lorentz algebras ${\cal U}^{(r,s)}_q(\SL)$.
In particular, for $r=1,~s=1$ its dual object is just the Gauss
form of the quantum Lorentz algebra \cite{WZ1}.

\bigskip
\noindent
{\bf 2. Quantum Lobachevsky space.}
The quantum Lobachevsky space  ${\bf L}_{\ka,q}$ is the associative $*$-algebra over
${\mathbb C}$ with a unity and three generators
$$
(z^*,H,z),~~H^*=H,~(z)^*=z^*, ~1^*=1.
$$
The commutation relations depend on two parameters $\ka,q\in{\mathbb C}$
\beq{3.9}
zH=\ka Hz,~~~~z^*H=\baka^{-1} Hz^*,
\eq
\beq{3.10}
zz^*=az^*z-bH^{-2}, ~~a=\left(\frac{\ka}{q}\right)^2,
~~b=\baka^{-1}\left(\frac{\ka}{q}\right)^2
\left(
1-q^2\frac{\baka}{\ka}\right)
\eq
The both parameters should be tune in a such a way that the
coefficients $a,b$ become real. It is happened in two cases\\
(i)$\ka$ and $q$ are real,\\
(ii) $|q|=1$ and ${\rm Arg}(q)={\rm Arg}(\ka)$.

\begin{predl}\label{p3.1}
${\bf L}_{\ka,q}$ is a right ${\cal U}^{(0,s)}_q(\SL)$-module\\
\end{predl}
{\sl Proof.}\\
Define the right action of ${\cal U}^{(0,s)}_q(\SL)$ on ${\bf L}_{\ka,q}$:
\beq{3.11}
\begin{array}{lll}
z^*.A=z^*, & H.A=q^{\oh}H, & z.A=q^{-1}z,\\
z^*.A^*=\left(\frac{\ka}{q}\right)^{\frac{2}{s}}z^*, &
H.A^*=\left(\frac{q}{\ka}\right)^{\frac{1}{s}} H, & z.A^*=z,\\
z^*.B=0, & H.B=0 & z.B=q^{-\oh} \\
z^*.C=q^{\frac{3}{2}}\ka^{-1}H^{-2}, & H.C=Hz, & z.C=-q^{\oh}z^2,
\end{array}
\eq
$$
x.B^*=0,  ~~~x.C^*=0~~{\rm for~any~}x\in{\bf L}_{\ka,q}.
$$
Similarly, the left action of $A^*,A,B^*,C^*$ takes the form
\beq{3.12}
\begin{array}{lll}
A^*.z^*=\bq^{-1}z^*, & A^*.H=\bq^\oh H, & A^*.z=z,\\
A.z^*=z^*, & A.H=\left(\frac{\bq}{\baka}\right)^{\frac1s} H,  &
A.z=\left(\frac{\baka}{\bq}\right)^{\frac{2}{s}}z\\
B^*.z^*=\bq^{-\oh}, & B^*.H=0, & B^*z=0\\
C^*.z^*=-\bq^{\oh}(z^*)^2, & C^*.H=z^*H, &
C^*.z=\bq^{\frac{3}{2}}\baka^{-1}H^{-2}.
\end{array}
\eq

The coproduct  in ${\cal U}^{(0,s)}_q(\SL)$ (\ref{3.4}),(\ref{3.5}) is compatible with
the commutation rules in ${\bf L}_{\ka,q}$ (\ref{3.9}), (\ref{3.10}).
Similar, the $*$-structures in ${\bf L}_{\ka,q}$ and in
${\cal U}^{(0,s)}_q(\SL)$ are also compatible.  It means that (\ref{3.11})
and (\ref{3.12}) are related as
$$
(x.K)^*=K^*.x^*.
$$
for $K=A,B,C$ and $x=z^*,H,z$. $\Box$
\bigskip

To eliminate the ambiguities related to the noncommutativity of ${\bf L}_{\ka,q}$
we consider only the ordered monomials putting $z^*$ on the left side, $z$ on
the right side and keeping $H$ in the middle of the monomials:
\beq{3.13}
:f(z^*,H,z):=\sum_{m,k,n}c_{m,k,n}w(m,k,n),~~
w(m,k,n)=z^{*m}H^kz^n.
\eq
The actions (\ref{3.11}),(\ref{3.12}) on $w(m,k,n)$
take the form
$$
w(m,k,n).A=q^{-n+\frac{k}{2}}w(m,k,n),
~~~~~~~w(m,k,n).A^*=\left(\frac{q}{\ka}\right)^{(-2m+k)/s}w(m,k,n),
$$
\beq{3.14}
w(m,k,n).B=q^{-n+\frac{k+1}{2}}
\frac{1-q^{2n}}{1-q^2}w(m,k,n-1),
\eq
$$
w(m,k,n).C=q^{n-\frac{3(k-1)}2}\ka^{k-1}
\frac{1-q^{2m}}{1-q^2}w(m-1,k-2,n)-
q^{-n+\frac{k+3}{2}}
\frac{1-q^{2n-2k}}{1-q^2}w(m,k,n+1),
$$
\bigskip
$$
A^*.w(m,k,n)=\bq^{-m+\frac{k}{2}}w(m,k,n),
~~~~~~~A.w(m,k,n)=\left(\frac{\bq}{\baka}\right)^{(-2n+k)/s}w(m,k,n),
$$
\beq{3.15}
B^*.w(m,k,n)=\bq^{-m+\frac{k+1}{2}}
\frac{1-\bq^{2m}}{1-\bq^2}w(m-1,k,n),
\eq
$$
C^*.w(m,k,n)=\bq^{m-\frac{3(k-1)}2}\baka^{k-1}
\frac{1-\bq^{2n}}{1-\bq^2}w(m,k-2,n-1)-
\bq^{-m+\frac{k+3}{2}}
\frac{1-\bq^{2m-2k}}{1-\bq^2}w(m+1,k,n).
$$
and
$$
w(m,k,n).\Om_q=q^{-k+1}\frac{(1-q^{k+1})^2}{(1-q^2)^2}w(m,k,n)+
$$
$$
+q^{-k+1}\ka^{k-1}\frac{(1-q^{2m})(1-q^{2n})}{(1-q^2)^2}w(m-1,k-2,n-1).
$$
Define the difference operators ${\cal D}_{z^*}, ~{\cal D}_z$ as
$$
{\cal D}_{z^*}f(z^*,H,z)=(z^*)^{-1}
\frac{f(z^*,H,z)-f(q^2z^*,H,z)}{1-q^2},
$$
$$
{\cal D}_zf(z^*,H,z)=\frac{f(z^*,H,z)-f(z^*,H,q^2z)}{1-q^2}z^{-1},
$$
and
$$
\bar{\cal D}_{z^*}f(z^*,H,z)=(z^*)^{-1}
\frac{f(z^*,H,z)-f(\bq^2z^*,H,z)}{1-\bq^2},
$$
$$
\bar{\cal D}_zf(z^*,H,z)=\frac{f(z^*,H,z)-f(z^*,H,\bq^2z)}{1-\bq^2}z^{-1}.
$$
Then the action (\ref{3.14}),(\ref{3.15}) of generators on the ordered functions
can be represent in the form
$$
f(z^*,H,z).A=f(z^*,q^\oh H,q^{-1}z),~~~~~~~~
f(z^*,H,z).A^*=
f(\left(\frac{\ka}{q}\right)^{\frac{2}{s}}z^*,
\left(\frac{\ka}{q}\right)^{\frac{1}{s}}H,z),
$$
\beq{3.16}
f(z^*,H,z).B=q^\oh {\cal D}_zf(z^*,q^\oh H,q^{-1}z),
\eq
$$
f(z^*,H,z).C=q^{\frac32}\ka^{-1}
{\cal D}_{z^*}f(z^*,q^{-\frac32}\ka H,q\ka^2z)H^{-2}+
\frac{q^{\frac32}}{1-q^2}
[f(z^*,q^{-\frac32}H,qz)-f(z^*,q^\oh H,qz)]z-
$$
$$
-q^{\frac32}{\cal D}_zf(z^*,q^\oh H,q^{-1}z)z^2,
$$
and
$$
A.f(z^*,H,z)=
f(z^*,\left(\frac{\baka}{\bq}\right)^{\frac1s}H,
\left(\frac{\baka}{\bq}\right)^{\frac{2}{s}}z),
~~~~~~~~A^*.f(z^*,H,z)=f(\bq^{-1}z^*,\bq^\oh H,z),
$$
\beq{3.17}
B^*.f(z^*,H,z)=\bq^\oh \bar{\cal D}_{z^*}f(\bq^{-1}z^*,\bq^\oh H,z),
\eq
$$
C^*.f(z^*,H,z)=
-\bq^{\frac32}(z^*)^2\bar{\cal D}_{z^*}f(\bq^{-1}z^*,\bq^\oh H,z)+
\frac{\bq^{\frac32}}{1-\bq^2}
z^*[f(\bq z^*,q^{-3/2}H,z)-f(\bq z^*,\bq^\oh H,z)]+
$$
$$
+\bq^{\frac32}\baka^{-1}H^{-2}\bar{\cal D}_zf(\bq\ka^2z^*,\bq^{-\frac32}\baka H,z).
$$

Note, that when $q\rar 1$ we come to (\ref{2.13})
$$
\p_q A\rar{\cal D}^3,~~~~B\rar{\cal D}^+,~~~~C\rar{\cal D}^-.
$$
\bigskip
The Haar functional $I_{{\bf L}_{\ka,q}} $ on ${\bf L}_{\ka,q}$ is a map
${\bf L}_{\ka,q}\rar{\mathbb C}$ that satisfies the following condition
$$
I_{{\bf L}_{\ka,q} }(f.u)=\ve (u)I_{{\bf L}_{\ka,q} }(f),~~f\in{\bf L}_{\ka,q},~~
u\in{\cal U}^{(0,s)}_q(\SL)
$$

\begin{predl}\label{p3.2}
 The Jackson integral for ordered $f(z^*,H,z)$ (\ref{3.13})
\beq{3.18}
I_{{\bf L}_{\ka,q} }(f)=\int\int\int :f(z^*,H,z)d_{q}z^*Hd_{q^\oh}Hd_{q}z:=
\eq
$$
=(1-q)^2(1-q^\oh\sum_{s,r,t\in {\bf Z}}
q^{s+ r+t}[:f(q^s,q^{\oh r},q^t):+:f(-q^s,q^{\oh r},-q^t):
$$
$$
+:f(q^s,q^{\oh r},-q^t):+:f(-q^s,q^{\oh r},q^t):].
$$
is the Haar functional on ${\bf L}_{\ka,q}$.
\end{predl}
{\sl Proof}\\
The main relation one should use to check (\ref{3.16}) is
$$
\int f(qx)d_qx=q^{-1}\int f(x)d_qx.
$$
Then the statement follows from the definitions of the actions (\ref{3.11}).
Moreover, from (\ref{3.16}), (\ref{3.17})) we find that the integral
(\ref{3.18}) is invariant under the action $f\rar u^*.f$. $\Box$
\vspace{10mm}

\section{Principle series of unitary representations}
\setcounter{equation}{0}

The principle series representations for quantum Lorentz groups were
considered in \cite{Pu,DDF}.
Here we present the formulae for principle series operators inspired
by the action of operators on the quantum Lobachevsky spaces. In
what follows we forget the coalgebraic structure and for this reason use
the notation ${\cal U}_q(\SL)$.

Consider the same space $V_{\nu,n}$ as in the classical case (\ref{2.3}). But now
we introduce another Hermitian form
\beq{4.1}
<f|g>=
q^{\frac{i\nu+n-1}{2}}
\bq^{\frac{i\nu-n+1}{2}}
\int\int \overline{
\left(
f(\bq^{\frac{i\nu-n+1}{2}}\bz, q^{\frac{i\nu+n-1}{2}}z)
\right)}
g(\bq^{\frac{i\nu-n+1}{2}}\bz,q^{\frac{i\nu+n-1}{2}}z)d\bz dz.
\eq

Define the action $\pi_\nu$ of ${\cal U}_q(\SL)$ on $V_{\nu,n}$. We assume that
${\cal U}_q({\rm SL}_2)$ acts only on the holomorphic part of $V_{\nu,n}$,
 while ${\cal U}^*_q({\rm SL}_2)$ acts on the
antiholomorphic part. Following our notations in Section 3 we write the holomorphic action
from the right hand side, while the antiholomorphic action from the left hand side.
The right action takes the form
$$
f(z).\pi_\nu(A)=q^{\frac{i\nu+n-1}{2}}f(q^{-1}z),
$$
\beq{4.2}
f(z).\pi_\nu (B)=q^{\frac{i\nu+n}{2}}
\frac{f(q^{-1}z)-f(qz)}{1-q^2}z^{-1},
\eq
$$
f(z).\pi_\nu (C)=\frac{z}{1-q^2}[-q^{\frac{i\nu+n+2}{2}}f(q^{-1}z)
+q^{-\frac{3(i\nu+n)}{2}+3}f(qz)],
$$
Similarly, for the left action we have
$$
\pi_\nu(A^*).f(\bz)=\bq^{\oh(i\nu-n-1)}f(\bq^{-1}\bz),
$$
\beq{4.3}
\pi_\nu (B^*).f(\bz)=\bq^{\frac{i\nu-n}{2}}
\frac{f(\bq^{-1}\bz)-f(\bq \bz)}{1-\bq^2}\bz^{-1},
\eq
$$
\pi_\nu (C^*).f(\bz)=\frac{\bz}{1-\bq^2}[-\bq^{\frac{i\nu-n+2}{2}}
f(\bq^{-1}\bz) +\bq^{-\frac{3(i\nu-n)}{2}+3}f(\bq \bz)].
$$

\begin{predl}\label{p4.1}
(i) The operators $\pi_\nu$ (\ref{4.2}),(\ref{4.3}) are irreducible
 representations of ${\cal U}_q(\SL)$;\\
(ii) for $\nu\in{\mathbb R}$ they are unitarity with respect to the form (\ref{4.1}):
$$
<g|f.\pi_\nu(u)>=
<\pi_\nu(S(u^*)).g|f>,
$$
$$
g,f\in V_{\nu,n},~~ u\in {\cal U}_q(\SL).
$$
\end{predl}
{\sl Proof.}
The direct calculations show that the operators $\pi_\nu(u)$ satisfy the same
commutation relations as $u$ (\ref{3.1}). Moreover, the representation
actions $\pi_\nu$ preserves the asymptotic (\ref{2.4}).

Another way to see that $\pi_\nu$ satisfy
the commutation relations for the case $n=0$ is to
consider the subspace $\ti{V}_\nu$ of ${\bf L}_{\ka,q}$ of the form
 $$
f(z^*,cH,z)=c^{i\nu-1}f(z^*,H,z).
$$
In other words, the basis functions in  $\ti{V}_\nu$ are
$w(m,i\nu-1,n)=z^{m*}H^{i\nu-1}z^n$.  Starting from (\ref{3.14}), (\ref{3.15})
 we define the action of ${\cal U}^{(0,s)}_q(\SL) )$ on $\ti{V}_\nu$ in a such way
that the right actions of $\pi_\nu(A,B,C)$ do not touch $z^*$-variable
($m=0$ in (\ref{3.14})). Similarly, the left acting operators  $\pi_\nu(A^*,B^*,C^*)$
 do not touch $z$-variable ($n=0$ in (\ref{3.15})).  Then $\ti{V}_\nu$ is
invariant with respect to $\pi_\nu$.
Since the dependence on $H$ is fixed one can put $H=1$. Then we come to
the representation (\ref{4.2}) and (\ref{4.3}), where $V_\nu$ is defined as
$V_\nu=\ti{V}_\nu|_{H=1}$.

Consider the action of the Casimirs (\ref{3.6}),(\ref{3.7})
 on $V_\nu$. They become the scalar operators
\beq{4.4}
\pi_\nu (\Om)=([-\frac{i\nu+n}{2}]_{q})^2{\rm Id},~
~\pi_\nu (\Om^*)=([-\frac{i\nu-n}{2}]_{\bq})^2{\rm Id}.
\eq
In this sense the representations $\pi_\nu$ are irreducible.

The unitarity is proved in Appendix. $\Box$
\bigskip

In the classical limit we come to the principle series representations of
$\SL$ since (\ref{2.7})
$$
\lim_{q\to 1}\pi_\nu(B)=T^+,~~
\lim_{q\to 1}\pi_\nu(C)=T^+,
$$
$$
\lim_{q\to 1} \pi_\nu(A)=(1+T^3),
$$
and the Hermitian form (\ref{4.1}) takes the classical form (\ref{2.6}).

It follows from (\ref{4.2}) and (\ref{4.3}) that $\nu$ is defined up to a constant.
Two representations
$(\nu,n)$ and $(\nu',n)$ are equivalent if
$$
\nu'=\nu+\frac{4\pi m}{\hbar},~~m\in{\mathbb Z}, ~~\hbar=\ln q,
$$
and we assume that
$$
\nu\in{\mathbb R},~~~0\leq\nu\le\frac{4\pi}{\hbar}.
$$

The case $n=0$ corresponds to the class one representations $V_{\nu,0}=V_\nu$.
It was demonstrated  in the
Proof of Proposition 4.1 that they are related to the quantum Lobachevsky space.
Moreover, as in the classical limit there exists a ${\cal U}_q(\SU)$ invariant state
in $ V_\nu$. It will be found in Section 6.
\vspace{10mm}

\section{The $q^2$-Bessel functions}
\setcounter{equation}{0}

We take here $q$ and $\ka$ to be real, and put $\ka=q^\de$.

Consider the eigenvalue problem for the Casimir (\ref{3.6}) acting on ${\bf L}_{\ka,q}$
$$
f_\nu^\de(z^*,H,z).\Om_q^\nu=([\frac{i\nu}{2}]_{q})^2f_\nu^\de(z^*,H,z).
$$
Thereby, $f_\nu^{(\de)}(z^*,H,z)\in V_\nu$.
It follows from (\ref{3.14})  that if
$$
\Om_q=\Om_q-q^{-i\nu+2}\frac{(1-q^{i\nu})^2}{(1-q^2)^2}
$$
then
$$
w(m,k,n).\Om_q^\nu=\frac{q^{-k+1}(1-q^{i\nu+k+1})(1-q^{-i\nu+k+1})}
{(1-q^2)^2}w(m,k,n)+
$$
\beq{5.1}
+\frac{q^{(k-1)(\de-1)}(1-q^{2m})(1-q^{2n})}{(1-q^2)^2}w(m-1,k-2,n-1).
\eq

We seek solutions of the equation
$$
f_\nu^{(\de)}(z^*,H,z).\Om_q^\nu=0,
$$
where $f_\nu^{(\de)}(z^*,H,z)$ takes the form
$$
f^{(\de)}(z^*,H,z)=e_{q^2}(i\mu(1-q^2)z^*)F_\nu^{(\de)}(H)e_{q^2}(i\mu(1-q^2)z),
$$
and $e_{q^2}$ is the $q^2$-exponential
$$
e_{q^2}(x)=\sum_{n=0}^\infty\frac{x^n}{(q^2,q^2)_n}=\frac1{(x,q^2)_\infty},
~~~~~~|x|<1.
$$
The $q^2$-exponentials $e_{q^2}(i\mu(1-q^2)z^*)$ and $e_{q^2}(i\mu(1-q^2)z)$
satisfy the equations
\beq{5.2}
{\cal D}_{z^*}e_{q^2}(i\mu(1-q^2)z^*)=i\mu e_{q^2}(i\mu(1-q^2)z^*),
~~{\cal D}_ze_{q^2}(i\mu(1-q^2)z)=i\mu e_{q^2}(i\mu(1-q^2)z).
\eq
On the other hand (\ref{5.1}) can be rewrite in the following
form
$$
e_{q^2}(i\mu(1-q^2)z^*)H^ke_{q^2}(i\mu(1-q^2)z).\Om_q^\nu=
$$
$$
=\frac{q^{-k+1}(1-q^{i\nu+k+1})(1-q^{-i\nu+k+1})}{(1-q^2)^2}
e_{q^2}(i\mu(1-q^2)z^*)F_\nu^{(\de)}(H)e_{q^2}(i\mu(1-q^2)z)+
$$
$$
+q^{(k-1)(\de-1)}
{\cal D}_{z^*}e_{q^2}(i\mu(1-q^2)z^*)H^{k-2}{\cal D}_ze_{q^2}(i\mu(1-q^2)z).
$$
It follows from the last expression and (\ref{5.2}) that $F_\nu^{(\de)}(H)$ is
the solution of the difference equation
\beq{5.3}
\frac{qF_\nu^{(\de)}(qH)-(q^{i\nu}+q^{-i\nu})F_\nu^{(\de)}(H)+
q^{-1}F_\nu^{(\de)}(q^{-1}H)}{(1-q^2)^2}-\mu^2q^{-\de-1}H^{-2}F_\nu^{(\de)}(q^{\de-1}H)=0.
\eq
Substitution $F_\nu^{(\de)}(H)=H^{-1}{\cal F}_\nu^{(\de)}(H^{-1})$ and $H^{-1}=x$ leads to
the difference equation for ${\cal F}_\nu^{(\de)}(x)$
\beq{5.4}
{\cal F}_\nu^{(\de)}(q^{-1}x)-(q^{i\nu}+q^{-i\nu}){\cal F}_\nu^{(\de)}(x)+
{\cal F}_\nu^{(\de)}(qx)-\mu^2q^{-2\de}(1-q^2)^2x^2{\cal F}_\nu^{(\de)}(q^{1-\de}x)=0
\eq
Define the second order difference operator
$$
\Delta_q \psi(x)=\frac{\psi(qx)-2\psi(x)+\psi(q^{-1}x)}{(q-q^{-1})^2}.
$$
Then (\ref{5.4}) takes the form
$$
\Delta_q{\cal F}_\nu^{(\de)}(x) - \mu^2q^{-2\de+2}x^2{\cal F}_\nu^{(\de)}(q^{1-\de}x)
=([\frac{i\nu}{2}]_{q})^2{\cal F}_\nu^{(\de)}(x).
$$
It is just the discrete analog of the quantum Liouville equation (\ref{2.24}) corresponding
to the  two-body case of the quantum open relativistic Toda lattices \cite{Ru} with
depending on the parameter $\de$ Hamiltonians
$$
\hat{H}=-\Delta_q+\mu^2q^{-2\de+2}x^2(T_q)^{1-\de},
$$
where $T_q$ is the shift operator $T_q\psi(x)=\psi(qx)$.
The wave-functions of these systems are among the solutions of (\ref{5.4}).
It follows from (\ref{5.4}) that the second ordered difference
equations arise only for three values $\de=0,1,2$. The other cases lead to  higher ordered
difference equations. In what follows we consider only these three cases.

Let
$$
{\cal F}_\nu^{(\de)}(x)=\sum_{k=0}^\infty\frac{a_k(1-q^2)^{2k}(\mu x)^{i\nu+2k}}
{(q^2,q^2)_k(q^{2i\nu+2},q^2)_k}
$$
Then from (\ref{5.4}) $a_k=(-1)^kq^{(2-\de)k(k+i\nu)-k\de}$, and
we find that the solutions are the modified $q^2$-Bessel functions
\cite{OR2}
$$
J_{i\nu}^{(j)}(2i\mu(1-q^2)q^{-\oh\de}x;q^2)=\frac{q^{-\oh i\de\nu}}
{\Gamma_{q^2}(i\nu+1)}
\sum_{k=0}^\infty\frac{(-1)^kq^{(2-\de)k(k+i\nu)-k\de}
(1-q^2)^{2k}(i\mu x)^{i\nu+2k}}{(q^2,q^2)_k(q^{2i\nu+2},q^2)_k}.
$$
Here
\beq{5.5}
j=\left\{
\begin{array}{lcl}
1 & {\rm for} & \de=2\\
2 & {\rm for} & \de=0\\
3 & {\rm for} & \de=1.\\
\end{array}
\right.
\eq
Three values of $j$ correspond to the known functions \cite{OR2,OR3}:
$$
I_{i\nu}^{(j)}(2\mu(1-q^2)q^{-\oh\de}x;q^2)=
e^{\oh\nu\pi}J_{i\nu}^{(j)}(2i\mu(1-q^2)q^{-\oh\de}x;q^2) -
$$
$$
- \left\{
\begin{array}{lc}
{\rm Jackson ~modified} ~q^2{\rm -Bessel~ function~of~type~1} \cite{Ja} &j=1\\
{\rm Hahn-Exton ~modified}~q^2{\rm -Bessel~
function ~\cite{Ko},\cite{I}}&j=3\\
{\rm Jackson ~modified} ~q^2{\rm -Bessel~function~of~type~2 ~\cite{Ja}}& j=2
\end{array}
\right.
$$

Obviously $I_{-i\nu}^{(j)}(2\mu(1-q^2)q^{-\oh\de}x;q^2)$ is the solution
of (\ref{5.4} also. The corresponding Wronskians are
$$
W^{(j)}(x;q^2)=I_{i\nu}^{(j)}(2\mu(1-q^2)q^{-\oh\de}x;q^2)
I_{-i\nu}^{(j)}(2\mu(1-q^2)q^{1-\oh\de}x;q^2)-
$$
$$
-I_{-i\nu}^{(j)}(2\mu(1-q^2)q^{-\oh\de}x;q^2)
I_{i\nu}^{(j)}(2\mu(1-q^2)q^{1-\oh\de}x;q^2)=
$$
$$
=\frac 1{\Gamma_{q^2}(i\nu+1)\Gamma_{q^2}(-i\nu+1)}\sum_{k=0}^\infty
\frac{q^{(2-\de)k(k-i\nu)-k\de}(1-q^2)^{2k}(\mu x)^{2k}}
{(q^2,q^2)_k}\times
$$
$$
\times\sum_{n=0}^k\left[
\begin{array}{c} k \\ n
\end{array}
\right]_{q^2}
\frac{q^{-2(2-\de)n(k-n-i\nu)}(q^{-i\nu+2k-2n}-q^{i\nu+2n})}
{(q^{2i\nu+2},q^2)_n(q^{-2i\nu+2},q^2)_n}.
$$
for $i\nu\ne n, ~~~n$ is an integer. Here
$$
\left[
\begin{array}{c} k \\ n
\end{array}
\right]_{q^2}=\frac{(q^2,q^2)_k}
{(q^2,q^2)_n(q^2,q^2)_{k-n}}.
$$
The interior sum $S^{(\de)}$ is equal
$$
S^{(\de)}=\left\{ \begin{array}{lc}
(-1)^kq^{-k(k+2i\nu-1)-i\nu}(1-q^{2i\nu}) & \de=0, ~~k\ge0,\\
q^{-i\nu}(1-q^{2i\nu}) & \de=1, k=0,\\
0 & \de=1, k>0,\\
-q^{-i\nu}(1-q^{2i\nu}) & \de=2,~~ k\ge0,
\end{array}
\right.
$$
So
$$
W^{(j)}=\left\{
\begin{array}{lc}
\frac{q^{-i\nu}(1-q^{2i\nu})}{\Gamma_{q^2}(i\nu+1)\Gamma_{q^2}(-i\nu+1)}
E_{q^2}(-4\mu^2(1-q^2)q^2x^2) & \de=0,\\
\frac{q^{-i\nu}(1-q^{2i\nu})}
{\Gamma_{q^2}(i\nu+1)\Gamma_{q^2}(-i\nu+1)} & \de=1,\\
\frac{q^{-i\nu}(1-q^{2i\nu})}{\Gamma_{q^2}(i\nu+1)\Gamma_{q^2}(-i\nu+1)}
e_{q^2}(4\mu^2(1-q^2)q^{-2}x^2) & \de=2,
\end{array}
\right.
$$
Thus, the modified $q^2$-Bessel functions
$I_{i\nu}^{(j)}$, $I_{-i\nu}^{(j)}$ form the fundamental system.
Though they are unbounded functions, and
our goal to find their linear combination providing the bounded solution as in
the classical case (\ref{2.22}).
\vspace{10mm}

\section{The q-Whittaker functions}
\setcounter{equation}{0}

In this Section we assume as before that $\ka=q^\de$ and $q$ are real.

Consider an analog of the  group element defined on
${\cal U}_q(\SL)\bigotimes{\bf L}_{\ka,q}$.
\footnote
{The group elements $g$ are defined in the tensor product
${\cal U}_q\otimes{\cal A}_q$  and satisfy the relation
$\D g\otimes id=g_{12}g_{23}$ \cite{Mi,MV}.}
\beq{6.1}
g(A|H)=e^{{\cal K}\otimes 2\phi},\qquad A=e^{\hbar{\cal K}},
\qquad H=e^\phi,\qquad \hbar=\ln q.
\eq
The action of $\pi_\nu(g(A|H))$ in $V_\nu$ (see (\ref{4.2}))
\beq{6.2}
\pi_\nu(g(A|H)): f(\bz,z)\rar  H^{i\nu-1}f(\bz,H^{-2}z).
\eq
We rewrite the Hermitian form (\ref{4.1}) as
\beq{6.3}
<f|g>=
q^{i\nu+n-1}
\int\int \overline{f(\bz,z)}g(q^{i\nu-n+1}\bz,q^{i\nu+n-1}z)d\bz dz.
\eq
By means of (\ref{6.3}) we
define the matrix element of $\pi_\nu(g(A,H))$ in  $V_{\nu}$
\beq{6.4}
F_{\nu}^{(\de)}(H)=<\psi_L|\psi_R.g(A|H)>=<S(g(A^*|H)).\psi_L|\psi_R>,
\eq
We drop here and in what follows the notion of the representation $\pi_\nu$.

Let $\psi_L,~\psi_R\in V_{\nu}$ satisfy the following conditions
\beq{6.5}
C^*.\psi_L(\bz,z)=\psi_L(\bz,z).B
\eq
\beq{6.6}
B^*.\psi_R.B=-\bar{\mu}\mu q^{1-\de}\psi_R.A^{2(\de-1)}.
\eq
In other words, from (\ref{4.2}), (\ref{4.3})
\beq{6.7}
B^*.\psi_R(\bz,z).B=
-\bar{\mu}\mu q^{(i\nu-2)(\de-1)}\psi_R(\bz,q^{2-2\de}z).
\eq
Thus $\psi_R$ can be considered as the {\sl quantum Whittaker vector}
for ${\cal U}_q(\SL)$.
\begin{predl}\label{p6.1}
$F_{\nu}^{(\de)}(H)$ satisfies (\ref{5.3})
\end{predl}
{\sl Proof.}
Since $F_{\nu}^{(\de)}(H)$ is the matrix element in the irreducible
space $V_{\nu}$ the Casimir operator $\Om_q$ (\ref{3.8}) acts as the
scalar (\ref{4.4})
\beq{6.8}
<\psi_L|\psi_R.g(A|H)(\frac{qA^2+q^{-1}A^{-2}-2}{(q-q^{-1})^2}+CB)>=
[-\frac{i\nu}{2}]_{q}^2<\psi_L|\psi_R.g(A|H)>.
\eq
It follows from (\ref{6.1}) that
\beq{6.9}
F_{\nu}^{(\de)}(qH)=<\psi_L|\psi_R.g(A|H)A^2>.
\eq
Then, the $A$-dependent part of left hand side in this
equation is boiled down to the second order difference operator
\beq{6.10}
F_{\nu}^{(\de)}(H).\left(\frac{qA^2+q^{-1}A^{-2}-2}{(q-q^{-1})^2}\right)=
\eq
$$
<\psi_L|\psi_R.g(A|H)
\left(
\frac{qA^2+q^{-1}A^{-2}-2}{(q-q^{-1})^2}
\right)
>=
\frac{qF_{\nu}^{(\de)}(qH)-2F_{\nu}^{(\de)}(H)+
q^{-1}F_{\nu}^{(\de)}(q^{-1}H)}{(q-q^{-1})^2}.
$$
On the other hand  (\ref{3.1}) provides the following relation
$$
g(A|H)C=CH^{-2}g(A|H).
$$
Thereby
$$
<\psi_L|\psi_R.g(A|H)CB>=H^{-2}<\psi_L|\psi_R.Cg(A|H)B>.
$$
>From   (\ref{6.5}) and the unitarity with respect to (\ref{6.3}) we have
$$
<\psi_L|\psi_R.g(A|H)CB>=H^{-2}<S(C^*).\psi_L|(\psi_R.B).g(A|H)>=
$$
$$
=H^{-2}<\psi_L.S(B)|(\psi_R.B).g(A|H)>=H^{-2}<\psi_L|(B^*.\psi_R.B).g(A|H)>.
$$
The condition (\ref{6.7}) and (\ref{6.9}) allow to rewrite the last expression as
$$
-\bar{\mu}\mu q^{1-\de}H^{-2}<\psi_L|(\psi_R.A^{2(\de-1)}).g(A|H)>=
-\bar{\mu}\mu q^{1-\de}H^{-2}<\psi_L|\psi_R.g(A|H)A^{2(\de-1)}>=
$$
$$
=-\bar{\mu}\mu q^{1-\de}H^{-2}<\psi_L|\psi_R.g(A|q^{\de-1}H)>=
-\bar{\mu}\mu q^{1-\de}H^{-2}F_\nu(q^{\de-1}H).
$$
Thus
$$
F_\nu^{(\de)}(H).CB=-\bar{\mu}\mu q^{1-\de}H^{-2}F_\nu(q^{\de-1}H).
$$
Substituting this equality with (\ref{6.10}) in (\ref{6.8}) we come to (\ref{5.3}).
$\Box$

\bigskip
\begin{cor}
 $F_\nu^{(\de)}(H)$ is the eigenfunction of the conjugate
Casimir operator $\Om_q^*$.
\end{cor}\label{c6.1}
{\sl Proof}. It follows from (\ref{4.4}) that
$$
\Om_q^*.F_\nu^{(\de)}(H)=([\frac{i\nu}{2}]_{\bq})^2F_\nu^{(\de)}(H).
$$
Since
$$
F_\nu^{(\de)}(H)=\left(F_\nu^{(\de)}(H)\right)^*=<g(A^*|H).\psi^*_R|\psi^*_L>
$$
one has
$$
\Om_q^*.F_\nu^{(\de)}(H)=
<(\frac{(qA^*)^2+q^{-1}(A^*)^{-2}-2}{(q-q^{-1})^2}+B^*C^*)g(A^*|H).\psi^*_R|\psi_L^*>.
$$
Then we can use the same conditions (\ref{6.5}), (\ref{6.6})
 because $(\psi.u)^*=u^*.\psi^*$. $\Box$

\bigskip

Now calculate (\ref{6.4}) explicitly.
\begin{predl}\label{p6.2}
The matrix element  $F_{\nu}^{(\de)}(H)=<\psi_L|g(A,H)|\psi_R>$ (\ref{6.4}) has
the form of the double integral
\beq{6.11}
F_{\nu}^{(\de)}(H)=q^{i\nu-1}
H^{i\nu-1}\int\int\frac{(-q^{2i\nu+4}\bar zz,q^2)_\infty}
{(-q^2\bar zz,q^2)_\infty}
\xi_1(q^{i\nu+1}z^*)\xi_2(q^{i\nu-1}H^{-2}z)dzdz^*,
\eq
where $\xi_1$, $\xi_2$ are defined
below (\ref{6.12}),\ref{6.13}).
\end{predl}
{\sl Proof.}
Assume that
$$
A^*.\psi_L(z^*,z).(A^{-1})=\psi_L(z^*,z).
$$
This condition along with (\ref{6.5}) means that $\psi_L$ is the
${\cal U}_q(\SU)$-invariant vector (see (\ref{3.3}))
in $V_\nu$.

It follows from (\ref{6.5}) that the state $\psi_L$ has the form
$$
\psi_L(\bz,z)=\sum_{m=0}^\infty\frac{b_m}{(q^2,q^2)_m}\bz^mz^m.
$$
Then,  from the second condition (\ref{6.5}) and from (\ref{4.2})
one can calculate the coefficients\\
$b_m=(-1)^mq^{2m}(q^{-2i\nu+2},q^2)_m$. Thereby
$$
\psi_L(\bz,z)=\sum_{m=0}^\infty\frac{(-1)^mq^{2m}(q^{-2i\nu+2},q^2)_m}
{(q^2,q^2)_m}(\bz z)^m=\frac{(-q^{-2i\nu+4}x,q^2)_\infty}
{(-q^2x,q^2)_\infty},
$$
where $x=|z|^2$.

Now consider the Whittaker state (\ref{6.6}),(\ref{6.7}).
Let us find it in the form of the product
$$
\psi_R(\bz,z)=\xi_1(\bz)\xi_2(z),
$$
$$
\xi_1(\bz)=\sum_{n=0}^\infty a_n\frac{(1-q^2)^n}{(q^2,q^2)_n}\bar{\mu}^n(\bz)^n,
~~~\xi_2(z)=\sum_{n=0}^\infty c_n\frac{(1-q^2)^n}{(q^2,q^2)_n}\mu^nz^n.
$$
Then from (\ref{6.7})
$$
a_nc_k=-\bar\mu\mu a_{n-1}c_{k-1}
q^{i\nu\de-2i\nu-2\de+1+(k-1)(3-2\de)+n}=
$$
$$
=a_0(i\mu)^nq^{\frac{n(n-1)}2+n(1-i\nu)}
c_0(i\mu)^kq^{\frac{k(k-1)}2(3-2\de)+k(1-i\nu)}q^{(i\nu\de-2\de+2)k},
$$
Then assuming $a_0=c_0=1$
\beq{6.12}
\xi_1(\bz)=\phantom._1\Phi_1(0;-q;q,-i\mu(1-q^2)q^{1-i\nu}\bz),
\eq
where $\phantom._r\Phi_s$ is the basic hypergeometric series \cite{GR}.
For $\xi_2(z)$ we have in the similar way
$$
c_k=(i\mu)^kq^{\frac{k(k-1)}2(3-2\de)+k(1-i\nu)+(i\nu\de-2\de+2)k}
$$
Thus we obtain
\beq{6.13}
\xi_2(z)=\left\{
    \begin{array}{lcl}
    \phantom._0\Phi_2(-;0,-q;q,-i\mu(1-q^2)q^{3-i\nu}z)&{\rm for}& \de=0\\
    \phantom._1\Phi_1(0;-q;q,-i\mu(1-q^2)qz)&{\rm for}& \de=1\\
    \phantom._3\Phi_1(0,0,0;-q;q,-i\mu(1-q^2)q^{i\nu-1}z)&{\rm for}& \de=2.\\
\end{array}
    \right.
\eq

Using (\ref{6.2}) and (\ref{6.3}) we define the matrix element $F_\nu(H)$ in
the form (\ref{6.11}). $\Box$

Let us show that the (\ref{6.11}) is the $q$-analog of the (\ref{2.22}).
If we rewrite the integral (\ref{6.11}) in the polar coordinates
$z=\rho e^{i\phi}$
$$
F_{\nu}^{(\de)}(H)=
2q^{i\nu-1}H^{i\nu-1}\int_0^\infty\frac{(-q^{2i\nu+4}\rho^2,q^2)_\infty}
{-q^2\rho^2,q^2)_\infty}\rho d\rho \int_{-\pi}^\pi\xi_1(q^{i\nu+1}\rho e^{-i\phi})
\xi_2(q^{i\nu-1}H^{-2}\rho e^{i\phi})d\phi,
$$
and represent the functions $\xi_1$ and $\xi_2$ as series, we obtain the
inner integral in the form
$$
\int_{-\pi}^\pi\sum_{n=0}^\infty\frac{q^{\frac{n(n-1)}2}(1-q^2)^n}{(q^2,q^2)_n}
(i\mu q^2\rho)^ne^{-in\phi}
\sum_{k=0}^\infty\frac{q^{\frac{3-2\de}2k(k-1)}(1-q^2)^k}{(q^2,q^2)_k}
(i\mu q^{i\nu\de-2\de+2}H^{-2}\rho)^ke^{ik\phi}d\phi.
$$
If $\de<\frac32$ these series converge uniformly with respect $\rho$ and
$\phi$, and we can integrate them term by term. We receive
$$
2\pi\sum_{n=0}^\infty\frac{(-1)^nq^{(2-\de)n^2}(1-q^2)^{2n}}{(q^2,q^2)_n^2}
(\mu q^{i\nu\de-\de+2}H^{-1}\rho)^{2n}=
2\pi J_0^{(j)}(2\mu(1-q^2)q^{(i\nu-1)\frac\de2+1}H^{-1}\rho;q^2).
$$
The $q^2$-Bessel function is the holomorphic one of $\rho$ if $\de<\frac32$
and it is the meromorphic function of $\rho$ with the ordinary poles
$\rho=\pm i\frac{q^{-i\nu}H}{2\mu(1-q^2)}$ if $\de=2$. So for $\de=0, ~1,
~2$
\beq{6.14}
F_\nu^{(\de)}(H)=4\pi q^{i\nu-1}H^{i\nu-1}\int_0^\infty
\frac{(-q^{2i\nu+4}\rho^2,q^2)_\infty}{(-q^2\rho^2,q^2)_\infty}
J_0^{(j)}(2\mu(1-q^2)q^{(i\nu-1)\frac\de2+1}H^{-1}\rho;q^2)\rho d\rho.
\eq
Let us prove now that this integral exists. Consider the first factor.
Using properties of the $q$-binomial formula \cite{OR3} we have
$$
\frac{(1+\rho^2)(-q^{2i\nu+4}\rho^2,q^2)_\infty}{(-q^2\rho^2,q^2)_\infty}=
\frac{(1+\rho^2)(-q^{2i\nu+2}\rho^2,q^2)_\infty}
{(1+q^{2i\nu+2}\rho^2)(-q^2\rho^2,q^2)_\infty}=
$$
$$
=\frac{(1+\rho^2)(q^{2i\nu},q^2)_\infty}{(1+q^{2i\nu+2}\rho^2)(q^2,q^2)_\infty}
\sum_{n=0}^\infty\frac{(q^{-2i\nu+2},q^2)_nq^{2i\nu n}}
{(q^2,q^2)_n(1+q^{2n+2}\rho^2)}.
$$
It means that this expression is absolutely integrable function. Consider the
second factor. It follows from \cite{OR2} that if $\rho\ne0$
$$
|\frac\rho{1+\rho^2}J_0^{(1)}(2\mu(1-q^2)q^{i\nu}H^{-1}\rho;q^2)|\le
$$
$$
\le\frac{A^{(1)}\sqrt\rho}{1+\rho^2}|e_q(-i\mu(1-q^2)q^{i\nu}H^{-1}\rho;q^2)+
ie_q(i\mu(1-q^2)q^{i\nu}H^{-1}\rho;q^2)|\le C^{(1)},
$$
$$
|\frac\rho{1+\rho^2}J_0^{(2)}(2\mu(1-q^2)qH^{-1}\rho;q^2)|\le
$$
$$
\le\frac{A^{(2)}\sqrt\rho}{1+q^2}|E_q(-i\mu(1-q^2)qH^{-1}\rho;q^2)+
iE_q(i\mu(1-q^2)qH^{-1}\rho;q^2)|\le C^{(2)}.
$$
For $j=3$ one has
$$
J_0^{(3)}(2x,q^2)=\frac{(qx^2,q^2)_\infty}{(q^2,q^2)_\infty}
\sum_{k=0}^\infty\frac{(-1)^kq^{k(k+1)}}{(q^2,q^2)_k(qx^2,q^2)_k}.
$$
So
$$
|\frac\rho{1+\rho^2}J_0^{(3)}(2\mu(1-q^2)q^{\frac{i\nu+1}2}H^{-1}\rho;q^2)|
\le C^{(3)},
$$
and (\ref{6.14}) converge absolutely.

Setting $q^{1+\frac{i\nu\de}2}\rho=r$ we obtain finally
$$
F_\nu^{(\de)}(H)=4\pi q^{i\nu(\de+1)}H^{i\nu-1}\int_0^{\infty e^{\frac{i\nu\de}2\ln q}}
\frac{(-q^{2i\nu+2-i\nu\de}r^2,q^2)_\infty}{(-q^{-i\nu\de}r^2,q^2)_\infty}
J_0^{(j)}(2\mu(1-q^2)q^{-\frac\de2}H^{-1}r;q^2)rdr=
$$

$$
=-\frac{8\pi q^{-\nu^2+\frac32i\nu\de}\ln q}
{(1-q^2)\G_{q^2}(i\nu+1)A_{i\nu}^{|1-\de|}}\mu^{i\nu}H^{-1}
K_{i\nu}^{(j)}(2\mu(1-q^2)q^{-\frac\de2}H^{-1};q^2),
$$
where
\beq{6.15}
A_{i\nu}=\sqrt{\frac{I_{i\nu}^{(2)}(2;q^2)}{I_{-i\nu}^{(2)}(2;q^2)}},
\eq
and $j=1,2,3$ are connected with $\de=2,0,1$ by relations (\ref{5.5})
(see \cite{R}). The curve of integration for $\de=2 ~(j=1)$ should satisfy the
condition $\nu\ln q\ne(k+\frac12)\pi$.
\vspace{10mm}

\section{The representation of the q-Bessel-Macdonald functions by the
Mellin-Barns integral}
\setcounter{equation}{0}

The representation (\ref{2.25}) is equivalent to the well-known Mellin transform
of the classical Bessel-Macdonald function
$$
\int_0^{i\infty}K_{i\nu}(2\mu H^{-1})(2\mu H^{-1})^{ip-1}dp=
-i2^{-ip-2}\G(\frac i2(p-\nu))\G(\frac i2(p+\nu)), ~~~
{\rm Im}\:s>{\rm Im}\:\nu\ge0,
$$
The inversion formula has the form
\beq{7.1}
K_{i\nu}(2\mu H^{-1}))=\frac1{8\pi}\int_{-\infty+i\sigma}^{\infty+i\sigma}
\G(\frac i2(p-\nu))\G(\frac i2(p+\nu))(2\mu H^{-1})^{-ip}dp.
\eq
Here we generalized these relations for the $q$ case (see also \cite{KhL})).
Consider the convergent integral
\beq{7.2}
G(x)=\frac1{2\pi}\int_{\sigma-i\infty}^{\sigma+i\infty}g(s)x^{-s}ds,
\eq
and assume that $G(x)$ satisfies (\ref{5.4}). Then
$$
\frac1{2\pi}\int_{\sigma-i\infty}^{\sigma+i\infty}q^{-s}(1-q^{s+i\nu})
(1-q^{s-i\nu})g(s)x^{-s}ds=
\frac{(1-q^2)^2}{2\pi}\int_{\sigma-i\infty}^{\sigma+i\infty}g(s)
q^{-2\de+(\de-1)s}x^{-(s-2)}ds,
$$
and we come to the recurrence relation for $g(s)$
$$
q^{-s}\frac{(1-q^{s+i\nu})(1-q^{s-i\nu})}{(1-q^2)^2}g(s)=
q^{(\de-1)s-2}g(s+2).
$$
Its solution has the form
\beq{7.3}
g(s)=q^{-\frac\de4s^2+\frac{2+\de}2s}\G_{q^2}(\frac{s+i\nu}2)\G_{q^2}(\frac{s-i\nu}2),
\eq
and for ${\rm Re}\:s>-{\rm Re}\:\nu\ge0$ the integral (\ref{7.2}) does
converge uniformly with respect to $x\in(1,\infty)$.

Since $I_{i\nu}^{(j)}$ and $K_{i\nu}^{(j)}$ form the fundamental system
of the solutions of (\ref{5.4}) \cite{R}
\beq{7.4}
G(x)=AI_{i\nu}^{(j)}(2\mu(1-q^2)q^{-\oh}x;q^2)+
BK_{i\nu}^{(j)}(2\mu(1-q^2)q^{-\oh}x;q^2).
\eq
Let $j=1$. It was shown in \cite{OR2} that
$I_{i\nu}^{(1)}(2\mu(1-q^2)q^{-\oh\de}x;q^2)$ is a meromorphic function with
the ordinary poles $x=\pm\frac{q^{-r+\oh\de}}{\mu(1-q^2)}, ~r=0,1,\ldots$,
while $K_{i\nu}^{(1)}(2\mu(1-q^2)q^{-\oh\de}x;q^2)$ and the left side of
(\ref{7.4}) are the holomorphic functions in the region ${\rm Re}\:x>0$.
So in this case $A=0$.

Let $j=2, 3$. Then, it follows from \cite{R1}, that $\lim_{x\to\infty}
I_{i\nu}^{(j)}(2\mu(1-q^2)q^{-\oh}x;q^2)=\infty,\\
\lim_{x\to\infty}K_{i\nu}^{(j)}(2\mu(1-q^2)q^{-\oh}x;q^2)=0$. Since the left hand side of
(\ref{7.4}) goes to zero if $x\to\infty$ we have for $j=2, 3 ~~A=0$.

It follows from (\ref{7.2}) that $g(s)$ (\ref{7.3}) is the inverse Mellin
transform of\\ $BK_{i\nu}^{(j)}(2\mu(1-q^2)q^{-\oh}x;q^2)$
\beq{7.5}
g(s)=B\int_0^\infty K_{i\nu}^{(j)}(2\mu(1-q^2)q^{-\oh}x;q^2)x^{s-1}dx.
\eq
Thus
$$
 K_{i\nu}^{(j)}(2\mu(1-q^2)q^{-\oh}x;q^2)=\frac1{2\pi B}
\int_{\sigma-i\infty}^{\sigma+i\infty}
q^{-\frac\de4s^2+\frac{2+\de}2s}\G_{q^2}(\frac{s+i\nu}2)\G_{q^2}(\frac{s-i\nu}2)
x^{-s}ds
$$
and we need to calculate $B$.
The $q$-Bessel-Macdonald functions satisfy the following relation \cite{R1}:
$$
\frac{q^{\oh\de}}{\mu(1+q)x}
\tilde{\cal D}_xx^{i\nu}K_{i\nu}^{(j)}(2\mu(1-q^2)q^{-\oh\de}x;q^2)=
-q^{-\frac{2-\de}2(i\nu-1)}x^{i\nu-1}
K_{\nu-1}^{(j)}(2(1-q^2)q^{1-\de}x;q^2),
$$
where
$$
\tilde{\cal D}_xf(x)=\frac{f(x)-f(qx)}{(1-q)x}.
$$
It means that
$$
\int_0^\infty K_{i\nu}^{(j)}(2\mu(1-q^2)q^{-\oh\de}x;q^2)x^{s-1}dx=
$$
\beq{7.6}
=-\frac{q^{\frac{2-\de}2(s+i\nu)}}{\mu(1+q)}
\int_0^\infty
\tilde{\cal D}_x\left(x^{i\nu+1}K_{i\nu+1}^{(j)}(2\mu(1-q^2)q^{-\oh\de}x;q^2)\right)
x^{s-i\nu-2}dx.
\eq
Let $f(x)$ be a smooth function and $f(x), f'(x)$ are integrable on $(0,\infty)$.
Then it can be proved that
$$
\int_0^\infty\tilde{\cal D}_xf(x)dx=\frac{\ln q}{1-q}f(0).
$$
Assume that $s=i\nu+2$. Remind that
\beq{7.7}
K_{i\nu}^{(j)}(2\mu(1-q^2)q^{-\oh\de}x;q^2)=
\eq
$$
=\frac12q^{\nu^2+i\nu}\G_{q^2}(i\nu)\G_{q^2}(1-i\nu)
\left[A_{i\nu}^{|1-\de|}I_{-i\nu}^{(j)}(2\mu(1-q^2)q^{-\oh\de}x;q^2)-
A_{-i\nu}^{|1-\de|}I_{i\nu}^{(j)}(2\mu(1-q^2)q^{-\oh\de}x;q^2)\right],
$$
where $A_{i\nu}$ is determined by (\ref{6.15}), and
\beq{7.8}
I_{i\nu}^{(j)}(2\mu(1-q^2)q^{-\oh\de}x;q^2)=\frac{q^{-\oh i\de\nu}}
{\Gamma_{q^2}(i\nu+1)}
\sum_{k=0}^\infty\frac{q^{(2-\de)k(k+i\nu)-k\de}
(1-q^2)^{2k}(\mu x)^{i\nu+2k}}{(q^2,q^2)_k(q^{2i\nu+2},q^2)_k}.
\eq
Then it follows from (\ref{7.5}) - (\ref{7.8})
$$
q^{\frac\de4\nu^2+i\frac{2-\de}2\nu+2}\G_{q^2}(i\nu+1)
=-B\frac{q^{(2-\de)(i\nu+1)}}{\mu(1+q)}
\int_0^\infty
\tilde{\cal D}_x\left(x^{i\nu+1}K_{i\nu+1}^{(j)}(2\mu(1-q^2)q^{-\oh\de}x;q^2)\right)dx=
$$
$$
=-\frac{B\ln q}{2(1-q^2)}\mu^{-i\nu-2}\G_{q^2}(i\nu+1)A_{i\nu+1}^{|1-\de|}
q^{(1+i\nu)(2-i\nu-\oh\de)}.
$$
Taking into account $A_{i\nu+1}=A_{i\nu}$ (see \cite{OR2}), we obtain
$$
B=-\frac{2(1-q^2)}{\ln q}\mu^{i\nu+2}A_{-i\nu}^{|1-\de|}q^{-\frac{4-\de}4\nu^2+\oh\de}.
$$
So we have proved
\begin{predl}\label{p7.1}
The $q$-Bessel-Macdonald functions can be represented by the Barns integral
$$
 K_{i\nu}^{(j)}(2\mu(1-q^2)q^{-\oh\de}x;q^2)=
$$
$$
=-\frac{\ln q\mu^{-i\nu-2}}{4\pi(1-q^2)}A_\nu^{|1-\de|}q^{\frac{4-\de}4\nu^2-\oh\de}
\int_{\sigma-i\infty}^{\sigma+i\infty}q^{-\frac\de4s^2+\frac{2+\de}2s}
\G_{q^2}(\frac{s+i\nu}2)\G_{q^2}(\frac{s-i\nu}2)x^{-s}ds,
$$
where $A_\nu$  is determined by (\ref{6.15}), and $j=1, 2, 3$ and
$\de=2, 0, 1$ are connected by (\ref{5.5}).
\end{predl}

Assuming $x=H^{-1}$ and changing the variable of integration $s=ip$ we obtain
\beq{7.9}
 K_{i\nu}^{(j)}(2\mu(1-q^2)q^{-\oh\de}H^{-1};q^2)=
\eq
$$
=-\frac{\ln q\mu^{-i\nu-2}}{4\pi i(1-q^2)}A_\nu^{|1-\de|}q^{\frac{4-\de}4\nu^2-\oh\de}
\int_{-\infty-i\sigma}^{\infty-i\sigma}q^{\frac\de4p^2+i\frac{2+\de}2p}
\G_{q^2}(\frac i2(p+\nu))\G_{q^2}(\frac i2(p-\nu))H^{ip}dp,
$$

Obviously the integral representation of $q^2$-Bessel-Macdonald function
(\ref{7.9}) is the $q$-analog of integral representation (\ref{7.1}).
The same representation for $|q|=1$ was obtained in \cite{KhL}.

\section*{Appendix}
{\sl Proof of Theorem }\\

Consider first the action of $\pi_\nu(A)$
$$
\int[f(z^*)]^*g(q^{i\nu-1}z)).\pi_\nu(A)d_{q^2}z=
\int[f(z^*)]^*q^{\frac{i\nu-1}2}g(q^{i\nu-1}q^{-1}z)d_{q^2}z=
$$
$$
=\int [q^{\frac{-i\nu+1}2}qf(qz)]^*g(q^{i\nu-1}z)d_{q^2}z=
\int [q^{-\frac{i\nu-1}2}f(qz^*)]^*g(q^{i\nu-1}z)d_{q^2}z=
$$
$$
=\int[\pi(S^{\cal F}A^*).(f(z^*)]^*g(q^{i\nu-1}z)d_{q^2}z.
$$
For $\pi_\nu(B)$ we have
$$
\int[(f(z^*))^*]^*g(q^{i\nu-1}z)).\pi_\nu(B)d_{q^2}z=
\int[f(z^*)]^*q^{\frac{i\nu}2}{\cal D}_zg(q^{-1}q^{i\nu-1}z)d_{q^2}z=
$$
$$
=-\int [q^{\frac{i\nu}{2}+1}{\cal D}_{z^*}f(q^{-1}z^*)]^*
g(q^{i\nu-1)}z)d_{q^2}z=
$$
$$
=\int[\pi_\nu S^{\cal F}(B^*).f(z^*)]^*g(q^{i\nu-1)}z)d_{q^2}z.
$$
Finally for $\pi_\nu (C)$
$$
\int[f(z^*)]^*g(q^{i\nu-1)}z).\pi_\nu (C)d_{q^2}z=
$$
$$
\int[f(z^*)]^*[-q^{\frac{3}{2}i\nu}g(q^{i\nu-2)}z)+
q^{-\oh i\nu+2}g(q^{i\nu)}z)d_{q^2}z=
$$
$$
=\int [-q^{-\frac{3}{2}i\nu+2}f(qz^*)z^*+
q^{\oh i\nu}f(q^{-1}z^*)z^*]^*g(q^{i\nu-1)}z)d_{q^2}z=
$$
$$
=\int [\pi_\nu S^{\cal F}(C^*).f(z^*)]^*g(q^{i\nu-1)}z)d_{q^2}z.
$$
\bigskip

\small{

\end{document}